\documentclass[12pt,a4paper]{article}
\usepackage{amsmath,amssymb,amscd}

\newtheorem{prop}{Proposition}

\newtheorem{theo}{Theorem}

\newtheorem{coro}{Corollary}
\newtheorem{lem}{Lemma}

\newtheorem{defi1}{Definition}
\newenvironment{defi}[1][]{\begin{defi1}[#1]\rm}{\end{defi1}}
\newenvironment{de}[1][]{\begin{defi1}\rm}{\end{defi1}}

\newtheorem{exe}{Example}

\newenvironment{ex}[1][]{\begin{exe}\rm}{\end{exe}}

\newtheorem{notati}{Notation}

\newenvironment{nota}[1][]{\begin{notati}\rm}{\end{notati}}

\newtheorem{rema}{Remark}

\newcommand{\B}{\ensuremath{\mathcal B}}

\newcommand{\C}{\ensuremath{\mathbb C}}

\newcommand{\N}{\ensuremath{\mathbb N}}

\newcommand{\Z}{\ensuremath{\mathbb Z}}
\newcommand{\Q}{\ensuremath{\mathbb Q}}

\newcommand{\s}{\ensuremath{\sigma}}

\newcommand{\pri}{\ensuremath{\smallsetminus}}

\begin{document}

\parindent=0cm
\hfill Communication in Algebra (to appear) 

\begin{center}
 {\LARGE\bf Coordinates of $R[x,y]$:\\
  Constructions and classifications.}\\
\vspace{.4cm}
{\large Eric Edo}\\
\vspace{.4cm}

\end{center}
\ \\
{\bf Abstract.} {\small Let $R$ be a PID. We construct and classify all coordinates of $R[x,y]$ of the form
$p_2y+Q_2(p_1x+Q_1(y))$ with $p_1,p_2\in{\rm qt}(R)$ and $Q_1,Q_2\in{\rm qt}(R)[y]$. From this construction
(with $R=K[z]$) we obtain non tame automorphisms $\s$ of $K[x,y,z]$ (where $K$ is a field of characteristic $0$) such that
the sub-group generated by $\s$ and the affine automorphisms contains all tame automorphisms.\\

}

\section{Introduction}

1) Let $K$ be a field. Due to the famous Jung-van der Kulk theorem (cf. [15] and [16]), the group of all automorphisms of the
$K$-algebra $K[x,y]$ is generated by the sub-groups of affine automorphisms and triangular automorphisms. Moreover, this group is the amalgamated
product of these two sub-groups along their intersection (cf. [12]). This result allows us both to construct
all automorphisms of $K[x,y]$ (by composition) and to classify them in terms of the length or the polydegree (see [13] and [14]).\\
2) There exist two classical ways to construct automorphisms of the $K$-algebra $K[x,y,z]$.
The first one is to compose affine automorphisms and triangular automorphisms (we obtain the so called tame automorphisms). Shestakov and Umirbaev
have proved that, if $K$ is a field of characteristic $0$, we do not obtain all automorphisms of $K[x,y,z]$ in this way (cf. [21]).
The second way consists to extend automorphisms of the $K[z]$-algebra $K[z][x,y]$ to obtain automorphisms of $K[x,y,z]$ fixing $z$ (called
$z$-automorphisms). This idea is developed in [8] (see also [17] \S 9.4).
We do not know whether all automorphisms of $K[x,y,z]$ may be obtained by composing $z$-automorphisms and affine ones.\\
3) In this context, it is natural to study the automorphisms of the $R$-algebra $R[x,y]$ thinking that $R$ is a PID, a UFD, a domain or
simply a general ring. When $R$ is a domain, an automorphism of $R[x,y]$ is roughly defined by one of his component (cf. Corollary~2).
This is the reason why we focus our attention on coordinates of $R[x,y]$.\\
4) In section~2, we introduce some classical notations and we recall well-known theorems: Nagata (see Theorem~\ref{dim1}), Russell-Sathaye
(see Theorem~\ref{RS}) and Shestakov-Umirbaev (see Theorem~\ref{SU}).\\
5) In section~3, we give the construction of automorphisms with one component of the form $d^{-1}\{q_2y+Q_2(q_1dx+Q_1(y))\}\in R[x,y]$
(cf. Theorem~\ref{RL2}).\\
6) In section~4, we develop the first elements of a theory of classification of these automorphisms. We distinguish which
are tame and which come from Russell-Sathaye construction (cf. Theorem~\ref{P5}). We prove that, if $R$ is a PID, all
polynomials of the form $p_2y+Q_2(p_1x+Q_1(y))\in R[x,y]$ with $p_1,p_2\in{\rm qt}(R)$ and $Q_1,Q_2\in{\rm qt}(R)[y]$ can
be written with the form considered in Theorem~\ref{RL2} (cf. Theorem~\ref{CL}).\\
7) There are many motivations to construct such automorphisms of $R[x,y]$:\\
- To construct non tame automorphisms of $K[z][x,y]$ (see for example [8])
which give non tame automorphisms of $K[x,y,z]$ using Shestakov-Umirbaev theorem (cf. [21]).\\
- To construct families of automorphisms of $\C[x,y]$ with generic length~3 to study the closure of the set automorphisms of $\C[x,y]$
with a fixed polydegree (see [10]).\\
- To give a criterion to check if there exists an automorphism of $R[x,y]$ sending $p_1x+Q_1(y)$ to $p_2x+Q_2(y)$ where
$p_1,p_2\in R\pri\{0\}$ and $Q_1,Q_2\in R[y]$ (see section~5). This question is linked with the work of Poloni (cf. [19])
about the classification of Danielewski hypersurfaces.\\
- To obtain non tame automorphisms $\s$ of $K[x,y,z]$ (where $K$ is a field of characteristic $0$) such that
the sub-group generated by $\s$ and the affine automorphisms contains all tame automorphisms (see Section~6).

\section{Preliminaries}

\begin{nota}
Let $R$ be a commutative ring.\\
a) We denote by $R^*$ the multiplicative group of units of $R$ and by $R^{\times}$ the set of non zero-divisors of $R$
(when $R$ is a domain, we have $R^{\times}=R\pri\{0\}$). We denote by ${\rm qt}(R)=(R^{\times})^{-1}R$ the total quotient ring of $R$
(when $R$ is a domain, ${\rm qt}(R)$ is the field of fractions of $R$).
We denote by $R^{\times}/R^*$ the quotient of $R^{\times}$ by the equivalence
relation $\sim$ defined by $r\sim s$ if and only if there exists $u\in R^*$ such that $r=us$, for all $r,s\in R^{\times}$.
We fix a subset ${\cal U}(R)$ of $R^{\times}$ such that, for all $r\in R^{\times}$, there exists a unique element $w_R(r)\in{\cal U}(R)$ such that
$r\sim w_{R}(r)$. For example, we can take for ${\cal U}(K[z])$ the set of unitary polynomials with $w_{K[z]}(P(z))={1\over {\rm lt}(P(z))}P(z)$
(where ${\rm lt}(P(z))$ is the leading term of $P(z)$), for all $P(z)\in K[z]^{\times}$ and we can take for ${\cal U}(\Z)$ the set $\N\pri\{0\}$
with $w_{\Z}(n)=|n|$, for all $n\in\Z^{\times}$. We denote by ${\rm Nil}(R)$ the ideal of nilpotent elements in $R$.\\
b) Let $f_1,\ldots,f_n\in R[x_1,\ldots,x_n]$ be polynomials, we denote by $(f_1,\ldots,f_n)$ the endomorphism $\s$ of the $R$-algebra
$R[x_1,\ldots,x_n]$ defined by $\s(x_i)=f_i$, for all $i\in\{1,\ldots,n\}$.\\
c) Let $f_1,f_2\in R[x,y]$, we denote by $\det(J\s)=(\partial_xf_1)(\partial_yf_2)-(\partial_xf_2)(\partial_yf_1)$
the Jacobian determinant of the endomorphism $(f_1,f_2)$.\\
d) We denote by ${\rm GA}_n(R)$ ($n\in\N\pri\{0\}$) the automorphisms group of the $R$-algebra $R[x_1,\ldots,x_n]$
(when $n=1$, $x_1=y$, when $n=2$, $(x_1,x_2)=(x,y)$ and when $n=3$, $(x_1,x_2,x_3)=(x,y,z)$).\\
e) We denote by $\pi=(y,x)\in{\rm GA}_2(R)$ or $\pi=(y,x,z)\in{\rm GA}_3(R)$ the automorphism exchanging $x$ and $y$.\\
f) We denote by ${\rm VA}_n(R)=\{F\in R[x_1,\ldots,x_n]\,;\,\exists\s\in{\rm GA}_n(R)\,;\,\s(x_n)=F\}$.\\
the set of $R$-{\it coordinates} (or $R$-{\it variables}) of $R[x_1,\ldots,x_n]$.\\
g) If $R$ is a domain, we use the following notations:\\
${\rm Aff}_n(R)=\{\s\in{\rm GA}_n(R)\,;\,\forall i\in\{1,\ldots,n\},\,{\rm deg}(\s(x_i))=1\}$
(for the affine automorphisms group),\\
${\rm BA}_n(R)=\{\s\in{\rm GA}_n(R)\,;\,\forall i\in\{1,\ldots,n\},\,\s(x_i)\in R^*x_i+R[x_{i+1},\ldots,x_n]\}$
(for the triangular automorphisms group) and\\
${\rm TA}_n(R)=<{\rm Aff}_n(R),{\rm BA}_n(R)>$ (for the tame automorphisms group).
\end{nota}

The following theorem is well-known and describes ${\rm VA}_1(R)$ (cf. [18]).

\begin{theo}[Nagata, 1972]\label{dim1} Let $P\in R[y]$ be a polynomial. The following assumptions are equivalent:\\
i) $P\in{\rm VA}_1(R)$,\\
ii) there exist $r\in R$, $u\in R^*$ and $N\in{\rm Nil}(R[y])$ such that:
$$P(y)=uy+r+N(y).$$
\end{theo}

{\bf Remark\ }\\
{\bf 1)} There exists an algorithm and even an explicit formula (see [12] Theorem~3.1.1 for the characteristic zero case and [1] Theorem~6.2
for the positive characteristic case) to compute the inverse of the automorphism
$\s\in{\rm GA}_1(R)$ defined by $\s(y)=P$. So we can say that ${\rm VA}_1(R)$ is well understood for every commutative ring $R$.\\
{\bf 2)} It is very important, in Theorem~\ref{dim1}, to consider a ring which is not a domain (since if $R$ is a domain ${\rm VA}_1(R)$ contains
only affine polynomials). Nevertheless, when we study ${\rm VA}_2(R)$, we often assume that $R$ is a domain, a UFD (unique factorization domain)
or even a PID (principal ideal domain) because the main applications are for $R=K[z]$, where $K$ is field.\\
{\bf 3)} If $R$ is a $\Q$-algebra, assumptions i) and ii) of Theorem~\ref{dim1} are equivalent to: iii) $P'(y)\in R[y]^*$.\\

The following corollaries of Theorem~\ref{dim1} are useful.

\begin{coro} We have: ${\rm VA}_2(R)+{\rm Nil}(R[x,y])={\rm VA}_2(R)$.
\end{coro}

\begin{coro}\label{Nagata2} Let $\s,\tau\in{\rm GA}_2(R)$ be automorphisms such that $\s(y)=\tau(y)$. We set: $Y=\s(y)=\tau(y)$.
We have: $\s(x)=u(Y)\tau(x)+r(Y)+N(x,y)$ with $r\in R[y]$, $u\in R[y]^*$ and $N\in{\rm Nil}(R[x,y])$.\\
If $R$ is a domain, we have: $\s(x)=u\tau(x)+r(Y)$ with $r\in R[y]$ and $u\in R^*$.
\end{coro}

{\bf Remark\ } This last corollary shows that in $R[x,y]$, where $R$ is a domain, a coordinate is exactly the orbit of an automorphism under the action
of the group of triangular automorphisms.\\

A first natural idea is to try to describe ${\rm VA}_2(R)$ using ${\rm VA}_1(R/I)$ for some (principal) ideals $I$ of $R$.

\begin{nota}
Let $I$ be an ideal of $R$. The canonical morphism $\phi_I\,:R\to R/I$ may be extended to a morphism from $R[x,y]$ to $(R/I)[x,y]$.
We still denote this morphism $\phi_I$. If $I=pR$ for some $p\in R$, we set: $\phi_p=\phi_{pR}$.
\end{nota}

The second natural idea is to define subclasses of $R[x,y]$ and to try to describe the intersection between ${\rm VA}_2(R)$
and each of these classes. The following definition come from [3]:

\begin{defi}[Berson, 2002]\label{Bc}
Let $p_1,\ldots,p_l\in R^{\times}$, and
$Q_1,\ldots,Q_l\in R[y]$. We define $F_l\in R[x,y]$ by induction on $l\in\N\pri\{0\}$:\\
1) $F_1(x,y)=p_1x+Q_1(y)$,\\
2) $F_2(x,y)=p_2y+Q_2(p_1x+Q_1(y))$,\\
3) $F_l(x,y)=p_lF_{l-2}+Q_l(F_{l-1})$, for all $l\ge 2$. We set (for $l\in\N\pri\{0\}$):
$$\B^l(R)=\{F_l(x,y)\,;\,p_1,\ldots,p_l\in R^{\times}\,,\,Q_1,\ldots,Q_l\in R[y]\}\subset R[x,y]$$
and $\B^0(R)=\{p_0y+q_0\,;\,p_0\in R^{\times},q_0\in R\}$.
\end{defi}

\begin{defi}[Rational classes]\label{Rc}
Let $l\in\N$ be an integer. We define:
$${\cal R}^l(R)={\cal B}^l({\rm qt}(R))\cap R[x,y].$$
\end{defi}
{\bf Remark\ } We assume that $R$ is a domain.\\
{\bf 1)} Let $l\in\N$ be an integer. A polynomial $F\in R[x,y]$ belongs to ${\cal R}^l(R)$ if and only if there exist
$\tau_1,\ldots,\tau_{l+1}\in {\rm BA}_2({\rm qt}(R))$ such that $F=\tau_1\pi\ldots\tau_l\pi\tau_{l+1}(y)$
and we can assume that $\tau_i(y)=y$ for $i\in\{1,\ldots,l\}$ (see [14]).\\
{\bf 2)} By Jung-van der Kulk theorem, we have ${\rm GA}_2({\rm qt}(R))={\rm TA}_2({\rm qt}(R))$. Using Bruhat decomposition in
${\rm Gl}_2({\rm qt}(R))$ we deduce that:
$${\rm VA}_2(R)=\bigcup_{l\in\N}{\cal R}^l(R)\cap{\rm VA}_2(R).$$
{\bf 3)} The first component of the Nagata automorphism (see point 5 in Remark after Theorem~\ref{RS}) is
in ${\cal R}^2(R)\cap{\rm VA}_2(R)$ but is not in $\bigcup_{l\in\N}{\cal B}^l(R)$
(see Proposition~\ref{N1}).\\
{\bf 4)} We have $\B^0(R)={\cal R}^0(R)$ and $\B^1(R)={\cal R}^1(R)$ but $\B^2(R)\varsubsetneq{\cal R}^2(R)$.\\

The description of $\B^1(R)\cap{\rm VA}_2(R)={\cal R}^1(R)\cap{\rm VA}_2(R)$ has been done by Russell and Sathaye (cf. [20]):

\begin{theo}[Russell, Sathaye, 1976]\label{RS}
Let $p_1\in R^{\times}$ be an non zero-divisor and let $Q_1\in R[y]$ be a polynomial. We set $F(x,y)=p_1x+Q_1(y)$.
The following assumptions are equivalent:\\
i) $F(x,y)\in{\rm VA}_2(R)$,\\
ii) $\phi_{p_1}(F(x,y))\in{\rm VA}_2(R/p_1R)$,\\
iii) $\phi_{p_1}(Q_1(y))\in{\rm VA}_1(R/p_1R)$.
\end{theo}

{\bf Remark\ }\\
{\bf 1)} Theorem~\ref{RS} is a particular case of Theorem~\ref{L2} ($Q_2(y)=y$).\\
{\bf 2)} Theorem~\ref{RS} is true for all $p_1\in R$ (including zero-divisors) as shown by Berson (see Theorem~1.2.6. in~[4]).\\
{\bf 3)} If we assume that $R$ is a domain and suppose iii). For all $u\in R^*$ and $Q_2\in R[y]$ such that $Q_2(Q_1(y))=y$ mod $p_1R[y]$.
We have:
$$\s=(\,(up_1)^{-1}(y-Q_2(F(x,y)))\,,\,F(x,y)\,)\in{\rm GA}_2(R)$$
and (by Corollary~2) every $\s\in{\rm GA}_2(R)$ such that $\s(y)=F(x,y)$ has this form.\\
{\bf 4)} With the notations of the previous point, we have: $\s\in{\rm TA}_2(R)$ if and only if there exist $a\in R^*$ and $b\in R$ such that
$Q_1(y)=ax+b$ modulo $p_1R[y]$ (see for example [9]).\\
{\bf 5)} The very classical example is the Nagata automorphism ($R=K[z]$ where $K$ is a field, $p_1=z^2$, $Q_1(y)=y+zy^2$, $Q_2(y)=y-zy^2$ and $u=-1$).\\

Let $K$ be a field. The groups ${\rm GA}_2(K[z])$ and $\{\s\in{\rm GA}_3(K)\,;\,\s(z)=z\}$
are canonically isomorphic (by the map $(f_1,f_2)\to (f_1,f_2,z)$ ). In this way, we can consider
${\rm TA}_2(K[z])$ as a sub-group of ${\rm GA}_3(K)$.

\begin{theo}[Shestakov, Umirbaev, 2004]\label{SU}\ \\
Let $K$ be a field of characteristic $0$.
$$\{\s\in{\rm GA}_3(K)\,;\,\s(z)=z\}\cap{\rm TA}_3(K)={\rm TA}_2(K[z]).$$
\end{theo}

This theorem (see [21]) is very strong because it's easy to check if a $z$-automorphisms is in ${\rm TA}_2(K[z])$
(see [9], see also [14] for an algorithm). In particular, we know, since [18], that the Nagata automorphism is not in ${\rm TA}_2(K[z])$
and Theorem~\ref{SU} implies is not in ${\rm TA}_3(K)$. An even stronger result is obtain in [22] where a conjecture from [18] is solved:

\begin{theo}[Umirbaev, Yu 2004]\label{UY}\ \\
Let $K$ be a field of characteristic $0$. Let $\s\in{\rm GA}_2(K[z])\pri{\rm TA}_2(K[z])$ be a wild automorphism of $K[z][x,y]$
then there exists no tame automorphism $\tau\in{\rm TA}_3(K)$ such that $\tau(y)=\s(y)$ (we say that $\s(y)$ is a wild coordinate
of $K[x,y,z]$).
\end{theo}

\section{Length~2 constructions.}

The description of $F(x,y)=p_2y+Q_2(p_1x+Q_2(y))\in\B^2(R)\cap{\rm VA}_2(R)$
begins in $[11]$ with the case $p_2=1$ and in $[9]$ for the case $p_1R+p_2R=R$.
The case $p_2=1$ also appear independently in [8] (see also [17] \S 9.4)
in the case $R=K[z]$ where $K$ is a field of characteristic $0$.
A compleat characterization is given in Theorem~\ref{L2}. In Theorem~\ref{RL2}
we study ${\cal R}^2(R)\cap{\rm VA}_2(R)$ (which strictly contains $\B^2(R)\cap{\rm VA}_2(R)$).
Theorem~\ref{L2} is a particular case of Theorem~\ref{RL2} ($d=1$, $p_1=q_1$ and $p_2=q_2$).

\begin{theo}\label{L2}
Let $p_1,p_2\in R^{\times}$ be non zero-divisors and let $Q_1,Q_2\in R[y]$ be polynomials.
We set: $F(x,y)=p_2y+Q_2(p_1x+Q_1(y))\in R[x,y]$. The following assumptions are equivalent:\\
i) $F(x,y)\in{\rm VA}_2(R)$,\\
ii) a) $\phi_{p_1}(F(x,y))\in{\rm VA}_2(R/p_1R)$ \hspace*{.1cm}and
\hspace*{.1cm} b) $\phi_{p_2}(F(x,y))\in{\rm VA}_2(R/p_2R)$,\\
iii) a) $\phi_{p_1}(p_2y+Q_2(Q_1(y)))\in{\rm VA}_1(R/p_1R)$ \hspace*{.1cm}and
\hspace*{.1cm} b) $\phi_{p_2}(Q_2(y))\in{\rm VA}_1(R/p_2R)$.
\end{theo}

\begin{ex}
We assume that $R=K[z]$, where $K$ is a field. We set: $p_1=z^2(z-1)$, $p_2=z$, $Q_1(y)=y+zy^2$ and $Q_2(y)=(z-1)(y+zy^2)$. The polynomial
$F(x,y)=p_2y+Q_2(p_1x+Q_1(y))$
is a coordinate by Theorem~\ref{L2}.
We have $p_1R+p_2R\ne R$. Moreover $F(x,y)$ is not a coordinate of length $"1+1"$ ({\it i. e.} is not a component of an automorphism composed by
two automorphisms constructed in Theorem~\ref{RS}, see Definition~5.3).
\end{ex}

\begin{theo}\label{RL2}
Let $d,q_1,q_2\in R^{\times}$ be non zero-divisors such that $dR+q_2R=R$ and let $Q_1,Q_2\in R[y]$ be polynomials such that
$\phi_d(q_2y+Q_2(Q_1(y)))=0$. We set: $F(x,y)=d^{-1}\{q_2y+Q_2(q_1dx+Q_1(y))\}\in R[x,y]$.
The following assumptions are equivalent:\\
i) $F(x,y)\in{\rm VA}_2(R)$,\\
ii) a) $\phi_{q_1}(F(x,y))\in{\rm VA}_2(R/q_1R)$ \hspace*{.2cm}and
\hspace*{.2cm} b) $\phi_{q_2}(F(x,y))\in{\rm VA}_2(R/q_2R)$.\\
iii) a) $\phi_{q_1}(F(0,y))\in{\rm VA}_1(R/q_1R)$ \hspace*{.2cm}and
\hspace*{.2cm} b) $\phi_{q_2}(Q_2(y))\in{\rm VA}_1(R/q_2R)$.
\end{theo}

Before proving Theorem~\ref{RL2}, we recall the following three classical lemmas. Lemma~\ref{ideal} is obvious,
Lemma~\ref{vaun} is a consequence of Lemma~1.11 in [3], and Lemma~\ref{vde} is a consequence of Lemma~1.1.8 in [12] (which is a corollary of the
formal inverse function theorem).

\begin{lem}\label{ideal}
Let $F\in R[x,y]$ and let $I$ be an ideal of $R$. If $F\in{\rm VA}_2(R)$ then $\phi_I(F)\in{\rm VA}_2(R/I)$.
\end{lem}

\begin{lem}\label{vaun} Let $Q\in R[y]$ and $H\in R[x,y]$. The following assumptions are equivalent:\\
i) $Q(H(x,y))\in{\rm VA}_2(R)$,\\
ii) $Q(y)\in{\rm VA}_1(R)$ and $H(x,y)\in{\rm VA}_2(R)$,\\
In particular, ${\rm VA}_2(R)\cap R[y]={\rm VA}_1(R)$.
\end{lem}

\begin{lem}\label{vde} Let $\s$ be an endomorphism of the $R$-algebra $R[x,y]$. We have $\s\in{\rm GA}_2(R)$ if and only if the following
two assumptions are fulfilled:\\
i) ${\rm det}(J\s)(0)\in R^*$,\\
ii) $\s\in{\rm GA}_2({\rm qt}(R))$.
\end{lem}

{\bf Proof (of Theorem~\ref{RL2}).}\\
{\it i)} $\Rightarrow$ {\it ii)}. This follows from Lemma~\ref{ideal}.\\
{\it ii)} $\Rightarrow$ {\it iii)}. In this part of the proof, we use Lemma~\ref{vaun}.\\
a) Since $\phi_{q_1}(F(x,y))\in{\rm VA}_2(R/q_1R)$, we have $\phi_{q_1}(F(0,y))\in{\rm VA}_2(R/q_1R)$ and
$\phi_{q_1}(F(0,y))\in{\rm VA}_1(R/q_1R)$ using Lemma~\ref{vaun}.\\
b) (This part of the proof is the only one where we use the hypothesis $dR+q_2R=R$).
Since $\phi_{q_2}(F(x,y))\in{\rm VA}_2(R/q_2R)$ and since $d$ is an invertible element modulo $q_2$,
we have $\phi_{q_2}(Q_2(q_1dx+Q_1(y)))\in{\rm VA}_1(R/q_2R)$ and
$\phi_{q_2}(Q_2(y))\in{\rm VA}_1(R/q_2R)$ using Lemma~\ref{vaun}.\\
{\it iii)} $\Rightarrow$ {\it i)}. In this part of the proof, we use Lemma~\ref{vde}.\\
By b), there exist $S,U\in R[y]$ such that $S(Q_2(y))=y+q_2U(y)$ $(1)$. There exists $V\in R[x,y]$ such that $S(q_2y+x)-S(x)=q_2V(x,y)$.\\
Changing $x$ to $Q_2(q_1dx+Q_1(y))$ in the previous equation, we have:\\
$S(dF(x,y))-S(Q_2(q_1dx+Q_1(y)))=q_2V(Q_2(q_1dx+Q_1(y)),y)$\\
\hspace*{5.9cm} $=q_2W(y)$ mod $q_1q_2 R[x,y]$ $(2)$,\\
where $W(y)=V(Q_2(Q_1(y)),y)\in R[y]$.\\
By a), there exists $T\in R[y]$ such that $T(F(0,y))=y$ mod $q_1R[y]$.
We have: $T(F(x,y))=T(F(0,y))=y$ mod $q_1R[y]$ $(3)$.\\
We set: $Q_3(y)=S(dy)-q_2\{U(Q_1(T(y)))+W(T(y))\}.$\\
Modulo $q_1q_2 R[x,y]$, we have:\\
$Q_3(F(x,y))=S(dF(x,y))-q_2\{U(Q_1(y))+W(y)\}$ \hspace*{2.5cm}(by $(3)$)\\
$\hspace*{2.2cm}=S(Q_2(q_1dx+Q_1(y)))-q_2U(Q_1(y))$ \hspace*{2.7cm}(by $(2)$),\\
$\hspace*{2.2cm}=q_1dx+Q_1(y)+q_2(U(q_1dx+Q_1(y))-U(Q_1(y)))$ \hspace*{.5cm}(by $(1)$),\\
$\hspace*{2.2cm}=q_1dx+Q_1(y)$.\\
Finally $q_1dx+Q_1(y)-Q_3(F(x,y))=0$ mod $q_1q_2 R[x,y]$ $(4)$.\\
We consider the following endomorphisms of ${\rm qt}(R)[x,y]$: $\tau_1=(q_1dx+Q_1(y),y)$,\\
$\tau_2=(d^{-1}\{q_2x+Q_2(y)\},y)$, $\tau_3=((q_1q_2)^{-1}(x-Q_3(y)),y)$, (we recall that $\pi=(y,x)$). We compute:
$$\s=\tau_1\pi\tau_2\pi\tau_3=(\,(q_1q_2)^{-1}\{q_1dx+Q_1(y)-Q_3(F(x,y))\}\,,\,F(x,y)\,).$$
By $(4)$, $\s$ is an endomorphism of $R[x,y]$. By the chain rule, we have:
${\rm det}(J\s)={\rm det}(J\tau_1){\rm det}(J\tau_2){\rm det}(J\tau_3)=q_1dd^{-1}q_2(q_1q_2)^{-1}=1$.
Using Lemma~\ref{vde}, we conclude that $\s\in{\rm GA}_2(R)$ and $F(x,y)\in{\rm VA}_2(R)$.\\

{\bf Remark\ } We use the notations of Theorem~\ref{RL2}.\\
{\bf 1)} We have: $\s^{-1}=\tau_3^{-1}\pi\tau_2^{-1}\pi\tau_1^{-1}$, hence:
$$\s^{-1}=(\,(q_1d)^{-1}\{q_1q_2(x+Q_3(y))-Q_1(G(x,y))\}\,,\,G(x,y)\,),$$
where $G(x,y)=q_2^{-1}dy-q_2^{-1}Q_2(q_1q_2x+Q_3(y))$. It is not easy to prove directly (without using Lemma~\ref{vde}) that the first component of $\s^{-1}$
is a polynomial in particular if $p_1R+p_2R\ne R$.\\
{\bf 2)} We set: $a=Q_2'(0)$ and $N(y)=Q_2(y)-ay$. We have: $Q_2(y)=ay+N(y)$.\\
If we do not assume $dR+q_2R=R$ (which is equivalent to $\phi_{q_2}(d)\in(R/q_2R)^*$)
but the weaker assumption $\phi_{q_2}(d)\in(R/q_2R)^{\times}$ then the condition ii)~b) in Theorem~\ref{RL2}
is equivalent to the following one:\\
{\it iii)'~b)} $\phi_{q_2}(N)\in{\rm Nil}(R/q_2R[y])$ and $\phi_{q_2R+aq_1R}(F(0,y))\in{\rm VA}_1(R/q_2R+aq_1R)$.\\
But we don't know whether {\it iii)~a)} and {\it iii)'~b)} imply {\it i)}.\\
We justify this equivalence.
Since $\phi_{q_2}(d)\in(R/q_2R)^{\times}$, we can consider the localization in $\phi_{q_2}(d)$ of the ring $R/q_2R$:\\
$R/q_2R[d^{-1}]=\{x\in{\rm qt}(R/q_2R)\,;\,\exists n\in\N\,,\,\phi_{q_2}(d)^nx\in R/q_2R\}$.\\
If $\phi_{q_2}(F(x,y))\in{\rm VA}_2(R/q_2R)$ then $Q_2(dq_1x+Q_1(y))\in{\rm VA}_2(R/q_2R[d^{-1}])$. By Lemma~\ref{vaun},
we deduce that $Q_2(y)\in{\rm VA}_1(R/q_2R[d^{-1}])$. This implies that $\phi_{q_2}(N)$ is nilpotent in
$R/q_2R[d^{-1}][y]$ and then in $R/q_2R[y]$.\\
Now, we assume that $\phi_{q_2}(N)\in{\rm Nil}(R/q_2R[y])$.\\
We have: $\phi_{q_2}(d^{-1}(N(Q_1(y))-N(dq_1x+Q_1(y)))\in{\rm Nil}(R/q_2R[x,y])$ and by Corollary~1:\\
$\phi_{q_2}(F(x,y))\in{\rm VA}_2(R/q_2R)$\\
$\Leftrightarrow$ $\phi_{q_2}(aq_1x+d^{-1}(q_2y+aQ_1(y)+N(Q_1(y))))\in{\rm VA}_2(R/q_2R)$\\
$\Leftrightarrow$ $\phi_{q_2R+aq_1R}(F(0,y))\in{\rm VA}_1(R/q_2R+aq_1R)$.\\
 The last $\Leftrightarrow$ is justified by Berson's improvement of Russell-Sathaye Theorem (see Remark~2 of Theorem~\ref{RS}).\\
{\bf 3)} If $R$ is a $\Q$-algebra, Theorem~\ref{RL2} is a consequence of the following deep result on locally nilpotent derivation
due to Daigle and Freudenburg (for the case $R$ a UFD, see [7]), Bhatwadekar and Dutta (for the case $R$ normal noetherian domain, see [6]),
Berson, van den Essen and Maubach (for the general case, see [5]):\\
Let $R$ be a $\Q$-algebra and let $F\in R[x,y]$ be a polynomial. We have: $F\in{\rm VA}_2(R)$ if and only if $F\in{\rm VA}_2({\rm qt}(R))$ and
$(\partial_xF)R+(\partial_yF)R=R[x,y]$.\\
Actually, we prove that {\it iii)} implies $(\partial_xF)R+(\partial_yF)R=R[x,y]$. We set: $I=(\partial_xF)R+(\partial_yF)R$,
$I_1=I+q_1R[x,y]$ and $I_2=I+q_3R[x,y]$. By iii)~a), using the remark 3) of Theorem~\ref{dim1}, we have: $\phi_{q_1}(F(0,y)')\in R/q_1R[y]^*$.
We deduce that $0=\phi_{I_1}(\partial_yF(x,y))=\phi_{I_1}(F(0,y)')\in R/I_1[x,y]^*$. Hence $I_1=R$ and $\phi_I(q_1)\in (R/I)^*~(1)$.
By iii)~b), using Remark~3) of Theorem~\ref{dim1}, we have: $\phi_{q_2}(Q_2'(y))\in R/q_2R[y]^*$. Using $(1)$ we obtain:
$0=\phi_{I_2}(\partial_xF(x,y))=\phi_{I_2}(q_1)\phi_{I_2}(Q_2'(q_1dx+Q_1(y)))\in R/I_2[x,y]^*$. Hence $I_2=R$ and $\phi_I(q_1)\in R/I^*~(2)$.
Finally, using $(1)$ and $(2)$, we have: $0=\phi_I(dq_1\partial_yF(x,y)-Q_1'(y)\partial_xF(x,y))=\phi_I(q_1q_2)\in R/I[x,y]^*$ and $I=R$.

\begin{ex} Let $R$ be a PID. We give a general family of examples. Let $d,p_1,p_2\in R^{\times}$ be such that
\hbox{${\rm gcd}(d,p_1)={\rm gcd}(d,p_2)={\rm gcd}(p_1,p_2)=1$}
and let $u,v\in R$ be such that $du+q_1v=1$. Let $Q_3,Q_4\in R[y]$ be two polynomials such that $\phi_{q_2}(Q_4(y))\in{\rm Nil}(R/q_2R)$.
We consider: $Q_1(y)=y+dQ_3(y)$,
$Q_2(y)=q_1\{(d-q_2v)y+dQ_4(y)\}$ and $F(x,y)=d^{-1}\{q_2y+Q_2(q_1dx+Q_1(y))\}$.
We have: $\phi_d(q_2y+Q_2(Q_1(y)))=\phi_d(q_2(1-q_1v)y)=0$. In one hand, we have:
$\phi_{q_1}(F(0,y))=\phi_{q_1}(d^{-1}q_2y)\in{\rm VA}_1(R/q_1R)$ and, on the other hand, we have:
$\phi_{q_2}(Q_2(y))=\phi_{q_2}(q_1dy+dQ_4(y))\in{\rm VA}_1(R/q_2R)$.
The assumption iii) of Theorem~\ref{RL2} is fulfilled and we deduce that $F(x,y)\in{\rm VA}_2(R)$. Let's now give explicit examples:\\
Consider the ring $K[z]$, where $K$ is a field, and take: $d=z^2$, $q_1=(z-1)^2$, $q_2=(z-2)^2$, $Q_1(y)=y+z^2y^2$ and\\
$Q_2(y)=(z-1)^2\{(-2z^3+8z^2-4z-4)y+z^2(z-2)y^2\}$.\\
Consider the ring $R=\Z$ of integers and take: $d=3$, $q_1=5$, $q_2=2$, $Q_1(y)=y+6y^2$ and $Q_2(y)=25y+30y^2$.
\end{ex}

\section{\bf Length 2 classification.}

In all this section, we assume that $R$ is a UFD.
\begin{de}
Let $p_1,p_2\in {\rm qt}(R)^{\times}$ and $Q_1,Q_2\in {\rm qt}(R)[y]$ be such that $F(x,y)=p_2y+Q_2(p_1x+Q_1(y))\in R[x,y]$.
By definition, we have $F(x,y)\in{\cal R}^2(R)$. If ${\rm deg}(Q_2)\le 0$, then $F(x,y)\in{\cal R}^0(R)$.
If ${\rm deg}(Q_2)=1$, then $F(x,y)\in{\cal R}^1(R)$. If ${\rm deg}(Q_1)\le 0$, then $\pi F(x,y)\in{\cal R}^1(R)$
(recall $\pi=(y,x)$). We say that $F\in R[x,y]$ is a \textit{rational length 2 polynomial}
if ${\rm deg}(Q_1)\ge 1$ and ${\rm deg}(Q_2)\ge 2$.
\end{de}

{\bf Remark\ } The first two differences between ${\cal R}^1(R)$ and ${\cal R}^2(R)$ are due to the following facts:
If $p_1,p_2,p_3,p_4\in {\rm qt}(R)^{\times}$ and $Q_1,Q_2,Q_3,Q_4\in {\rm qt}(R)[y]$, we have:\\
\textbf{1)} $p_1x+Q_1(y)=p_2x+Q_2(y) \Leftrightarrow p_1=p_2$ and $Q_1=Q_2$ (the parameters of a polynomial in ${\cal R}^1(R)$ are unique).\\
\textbf{2)} $p_1x+Q_1(y)\in R[x,y]\Leftrightarrow p_1\in R$ and $Q_1\in R[y]$ (${\cal R}^1(R)={\cal B}^1(R)$).\\
\textbf{3)} $p_2y+Q_2(p_1x+Q_1(y))=p_4y+Q_4(p_3x+Q_3(y))\not\Rightarrow p_1=p_3,p_2=p_4$ and $Q_1=Q_3,Q_2=Q_4$ (the parameters of a rational length 2 polynomial
are not unique).\\
\textbf{4)} $p_2y+Q_2(p_1x+Q_1(y))\in R[x,y]\not\Rightarrow p_1,p_2\in R$ and $Q_1,Q_2\in R[y]$ (the parameters of a rational length 2 polynomial
are not always in the ring $R$).\\

The following proposition shows that there exists rational length 2 polynomial which are coordinate but are not in the Berson classes.

\begin{prop}\label{N1}
Let $K$ be a field of characteristic zero. We consider $p_1=z^2$ and $p_2=-z^{-2}$ in ${\rm qt}(K[z])=K(z)$,
$Q_1(y)=y+zy^2$ and $Q_2(y)=z^{-2}(y-zy^2)$ in ${\rm qt}(K[z])[y]=K(z)[y]$. The following properties hold:\\
1) $N(x,y)=p_2y+Q_2(p_1x+Q_1(y))=x-2y(zx+y^2)-z(zx+y^2)^2$ is a rational length 2 polynomial,\\
2) $(N(x,y),p_1x+Q_1(y))\in{\rm GA}_2(K[z])$ (this the Nagata automorphism),\\
3) For all $l\in\N$, we have $N(x,y)\not\in{\cal B}^l(K[z])$.
\end{prop}

{\bf Proof.} 1) is trivial and 2) is classical (see [18]). For a proof of 3) see [4]~Proposition~2.1.15.\\

The aim of the following definitions is to give some canonical parameters for a rational length 2 polynomial.

\begin{de}\ \\
a) We denote by $L_2(R)$ the set of all quadruplets $(p_1,p_2,Q_1,Q_2)$ where $p_1,p_2\in{\rm qt}(R)^{\times}$ and $Q_1,Q_2\in{\rm qt}(R)[y]$ are such that
${\rm deg}(Q_1)\ge 1$, \hbox{${\rm deg}(Q_2)\ge 2$} and $p_2y+Q_2(p_1x+Q_1(y))\in R[x,y]$.\\
b) We define an equivalence relation $\simeq$ between quadruplets in $L_2(R)$ in the following way:
$(p_1,p_2,Q_1,Q_2)\simeq(p_3,p_4,Q_3,Q_4)$ if there exists $r\in R$ such that $p_2y+Q_2(p_1x+Q_1(y))+r=p_4y+Q_4(p_3x+Q_3(y))$.\\
c) A quadruplet $(p_1,p_2,Q_1,Q_2)\in L_2(R)$ is said to be \textit{reduced} if the following conditions hold:
$Q_1(0)=Q_2(0)=0$, $p_1\in {\cal U}(R)$, $Q_1(y)\in R[y]$ and $\gcd(p_1,Q_1(y))=1$. We denote by $L_2^{\rm red}(R)$ the subset of all reduced quadruplets.
\end{de}

\begin{prop}\label{P0}
Every quadruplet in $L_2(R)$ is equivalent to a unique reduced quadruplet.
\end{prop}

{\bf Proof.} Let $(p_1,p_2,Q_1,Q_2)\in L_2(R)$ be a quadruplet.\\
1) We change $Q_1(y)$ to $Q_1(y)-Q_1(0)$ and $Q_2(y)$ to $Q_2(y+Q_1(0))$.\\
2) We change $Q_2(y)$ to $Q_2(y)-Q_2(0)$.\\
3) Let $m\in R^{\times}$ be the smallest common multiple of the denominators of $p_1$ and all the coefficients of $Q_1(y)$.
We change $Q_2(y)$ to $Q_2({y\over m})$, $p_1$ to $mp_1$ and $Q_1(y)$ to $mQ_1(y)$.\\
4) Let $u\in R^*$ be such that $p_1=u\,w_R(p_1)$ ($w_R(p_1)\in{\cal U}(R)$ see Notation~1~a). We change $Q_2(y)$ to $Q_2(uy)$, $p_1$ to $w_R(p_1)$
and $Q_1(y)$ to $u^{-1}Q_1(y)$.\\
After these $4$ modifications, we obtain a reduced quadruplet of $L_2^{\rm red}(R)$ which is equivalent to $(p_1,p_2,Q_1,Q_2)$.\\
Now, let $(p_1,p_2,Q_1,Q_2),(p_3,p_4,Q_3,Q_4)\in L_2^{red}(R)$ be equivalent reduced quadruplets. There exists $r\in R$ such that
$p_2y+Q_2(p_1x+Q_1(y))+r=p_4y+Q_4(p_3x+Q_3(y))$. Taking $x=y=0$, we obtain $r=0$. After changing $x$ to $p_1^{-1}x-p_1^{-1}Q_1(y)$ we have:\\
$(*)\,\,p_2y+Q_2(x)=p_4y+Q_4(p_3p_1^{-1}x+Q_3(y)-p_3p_1^{-1}Q_1(y))$.\\
Taking $x=0$ in $(*)$, we have: $p_2y=p_4y+Q_4(Q_3(y)-p_3p_1^{-1}Q_1(y))$. Since ${\rm deg}(Q_4)\ge 2$ and $Q_1(0)=Q_3(0)=0$, we deduce
$p_2=p_4$ and $p_1Q_3(y)=p_3Q_1(y)$. Since ${\rm gcd}(p_1,Q_1(y))={\rm gcd}(p_3,Q_3(y))=1$ this implies $p_1\sim p_3$ and $p_1=p_3$ (because
$p_1,p_2\in{\cal U}(R)$) and $Q_1(y)=Q_3(y)$. Now $(*)$ gives $Q_2(y)=Q_4(y)$.

\begin{prop}\label{P1}
Let $(p_1,p_2,Q_1,Q_2)\in L_2^{\rm red}(R)$ be a reduced quadruplet. We set $F(x,y)=p_2y+Q_2(p_1x+Q_1(y))$.
Then following properties hold:\\
1) $F(0,0)=0$,\\
2) $p_1p_2\in R$,\\
3) $p_2\in R$ if and only if $Q_2(y)\in R[y]$,\\
4) There exist $d\in {\cal U}(R)$, $q_1,q_2\in R^{\times}$ and a polynomial $\tilde{Q}_2\in R[y]$ unique such that
${\rm gcd}(d,q_2)=1$, $p_1=dq_1$, $p_2=d^{-1}q_2$ and $Q_2(y)=d^{-1}\tilde{Q}_2(y)$.
\end{prop}

{\bf Proof.}
1) $F(0,0)=Q_2(Q_1(0))=0$.\\
2) We have: $p_1p_2=p_1\partial_yF(x,y)-Q_1'(y)\partial_xF(x,y)\in R[x,y]\cap {\rm qt}(R)=R$.\\
3) If $Q_2(y)\in R[y]$ then $p_2y=F(x,y)-Q_2(p_1x+Q_1(y))\in R[x,y]\cap{\rm qt}(R)[y]=R[y]$ and $p_2\in R$.
Conversely, if $p_2\in R$ then $Q_2(p_1x+Q_1(y))\in R[x,y]$. By contradiction, let us assume $Q_2(y)\not\in R[y]$.
Let $m\in R^{\times}\pri R^*$ be the smallest common multiple of all denominators of coefficients in $Q_2(y)$.
Then $(mQ_2)(p_1x+Q_1(y))\in mR[x,y]$ and $mQ_2\in R[y]$ with ${\rm gcd}(mQ_2(y))=1$ and $p_1x+Q_1(y)\in R[x,y]$ with ${\rm gcd}(p_1x+Q_1(y))=1$
which is impossible.\\
4) There exist $d\in{\cal U}(R)$ and $q_2\in R^{\times}$ such that $p_2=d^{-1}q_2$ and ${\rm gcd}(d,q_2)=1$. By 2),
we have: $p_1p_2\in R$, hence $p_1q_2\in dR$ and this implies $p_1\in dR$ since ${\rm gcd}(d,q_2)=1$. In other words, there exists $q_1\in R^{\times}$
such that $p_1=dq_1$.\\
We have: $q_2y+(dQ_2)(p_1x+Q_1(y))=d(p_2y+Q_2(p_1x+Q_1(y)))\in R[x,y]$, hence $(p_1,q_2,Q_1,dQ_2)\in L_2^{\rm red}(R)$.
By 3), (since $q_2\in R$), we have: $\tilde{Q}_2(y)=dQ_2(y)\in R[y]$.

\begin{theo}\label{CL}
 Assume $R$ is PID. Then in Theorem~\ref{RL2}, we have constructed all coordinates of rational length~2.
\end{theo}

{\bf Proof.} This result follows from 4) of Proposition~\ref{P1} and Theorem~\ref{RL2} and the fact
that ${\rm gcd}(d,q_2)=1$ is equivalent to $dR+q_2R=R$ when $R$ is a PID.\\

For all the remaining of this section, we fix $(p_1,p_2,Q_1,Q_2)\in L_2^{\rm red}(R)$ a reduced quadruplet such that
$F(x,y)=p_2y+Q_2(p_1x+Q_1(y))\in{\rm VA}_2(R)$.

\begin{de}\ \\
1) We say that $F$ \emph{is tame} if there exists a tame automorphism $\s$ of $R[x,y]$ such that $\s(y)=F$.\\
2) We say that $F$ \emph{has a mate of length~$1$} if there exists $G\in{\cal B}^1(R)$ such that $(F,G)\in{\rm GA}_2(R)$.\\
3) We say that $F$ \emph{has length $"1+1"$} if there exists $\s,\tau\in{\rm GA}_2(R)$ such that $\s(y),\tau(y)\in{\cal B}^1(R)$ and $\s\tau(y)=F$.\\
\end{de}

{\bf Remark\ } Coordinate of length $"1+1"$ may be constructed with the help of Theorem~\ref{RS} and may be considered as
trivial from the length~2 point of view.

\begin{theo}\label{P5}
1) $F$ is tame if and only if there exist $\s,\tau\in{\rm GA}_2(R)$ such that $\s(x),\s(y),\tau(x),\tau(y)\in{\cal B}^1(R)$ and $\s\tau(y)=F$.\\
2) $F$ has a mate of length~$1$ if and only if $p_1p_2\in R^*$.\\
3) $F$ has length $"1+1"$ if and only if $\phi_{p_1}(Q_1(y))\in{\rm VA}_1(R/q_1R)$.\\
4) if $F$ is tame or if $F$ has a mate of length~$1$ then $F$ has length $"1+1"$.
\end{theo}

{\bf Proof.}\\
1) We assume that $F$ is tame. Let $\rho$ be a tame automorphism of $R[x,y]$ such that $\rho(y)=F$. We set: $\tau_1=(p_1x+Q_1(y))$ and $\tau_2=(p_2x+Q_2(y))$,
$\tau_1$ and $\tau_2$ are triangular automorphisms
of ${\rm qt}(R)[x,y]$. By Corollary~2, there exists $\tau_3$ a triangular automorphism of ${\rm qt}(R)[x,y]$ such that $\rho=\tau_1\pi\tau_2\pi\tau_3$.
The amalgamated structure of ${\rm GA}_2({\rm qt}(R))$ implies that $\rho=b_1a_1b_2a_2b_3$ where $a_i$ (resp. $b_i$) are affine (resp. triangular)
automorphisms of $R[x,y]$. Let $b_3'$ be such that $b_3'(x)=x$ and $b_3'(y)=b_3$ ($b_3'$ is an affine automorphism). We set: $\s=b_1a_2$ and
$\tau=b_2a_2b_3'$. We have: $\s\tau(y)=b_1a_1b_2a_2b_3'(y)=b_1a_1b_2a_2b_3(y)=\rho(y)=F$ and we easily verify that
$\s(x),\s(y),\tau(x),\tau(y)\in{\cal B}^1(R)$. Conversely, if there exist $\s,\tau\in{\rm GA}_2(R)$ such that
$\s(x),\s(y),\tau(x),\tau(y)\in{\cal B}^1(R)$ and $\s\tau(y)=F$ then $\s$ and $\tau$ are tame automorphisms and $F$ is tame.\\
2) At first, we assume that there exist $p_3\in R^{\times}$ and $Q_3\in R[y]$ such that $\s=(F(x,y),p_3x+Q_3(y))\in{\rm GA}_2(R)$.
Composing $\s$ with a translation we can assume that $Q_3(0)=0$.
We consider $\tau=(p_3x+Q_3(y),y)\in{\rm GA}_2({\rm qt}(R))$. We have:
$\pi\tau^{-1}\s=(p_2x+Q_2(p_1p_3^{-1}y+Q_1(x)-p_1p_3^{-1}Q_3(x)),y)\in{\rm GA}_2({\rm qt}(R))$.
Corollary~\ref{Nagata2} gives $p_2x+Q_2(p_1p_3^{-1}y+Q_1(x)-p_1p_3^{-1}Q_3(x))\in {\rm qt}(R)x+{\rm qt}(R)[y]$. Since ${\rm deg}(Q_2)\ge 2$ and $Q_1(0)=Q_3(0)=0$,
we deduce $p_1Q_3(y)=p_3Q_1(y)\ (*)$ and then $p_1Q_3'(y)=p_3Q_1'(y)\ (**)$.\\
Since ${\rm gcd}(p_1,Q_1(y))=1$, $(*)$ implies that there exists $u\in R$ such that $p_3=up_1$. Since $u(p_1x+Q_1(y))=p_3x+Q_3(y)=\s(y)\in{\rm VA}_2(R)$,
we deduce that $u\in R^*$ (when $R$ is a domain all coordinates are irreducible polynomials). Using $(**)$ we obtain:\\
${\rm det}(J\s)=Q_3'(y)p_1Q_2'(p_1x+Q_1(y))-p_3(p_2+Q_1'(y)Q_2'(p_1x+Q_1(y)))$\\
\hspace*{1.5cm}$=-p_2p_3\in R^*$.\\
Finally $p_2p_3=up_1p_2\in R^*$ implies $p_1p_2\in R^*$.\\
Conversely, we assume that $u=p_1p_2\in R^*$. Changing $p_1$ to $u^{-1}p_1$, $Q_1(y)$ to $u^{-1}Q_1(y)$ and $Q_2(y)$ to $Q_2(uy)$, one can assume $u=1$.\\
We set: $Q_3(y)=-p_1Q_2(y)\in {\rm qt}(R)[y]$. We have: $Q_3(p_1x+Q_1(y))=y-p_1F(x,y)\in R[x,y]$. Since ${\rm gcd}(p_1,Q_1(y))=1$, this implies
$Q_3(y)\in R[y]$. From $y-Q_3(p_1x+Q_1(y))=p_1F(x,y)\in p_1R[x,y]$, we deduce that $\phi_{p_1}(Q_1)$ and $\phi_{p_1}(Q_3)$ are inverse
in ${\rm GA}_1(R/p_1R)$ and $(F(x,y), p_1x+Q_1(y))\in{\rm GA}_2(R)$ by Theorem~\ref{RS}.\\
3) We assume that there exist $\s,\tau\in{\rm GA}_2(R)$ such that $\s(y),\tau(y)\in{\cal B}^1(R)$ and $\s\tau(y)=F$.
There exist $p_3,p_4\in R^{\times}$ and $Q_3(y),Q_4(y)\in R[y]$ such that $\s(y)=p_3x+Q_3(y)$ and $\tau(y)=p_4x+Q_4(y)$.
Let $u\in R^*$ be such that $p_3=uw_R(p_3)$. Changing $(\s,\tau)$ to $(\s\rho,\rho^{-1}\tau)$ where $\rho=(x,u(y-Q_3(0)))$ we can assume
that $p_3\in{\cal U}(R)$ and $Q_3(0)$. Since $\s(y)\in{\rm VA}_2(R)$, Theorem~\ref{RS} implies ${\rm gcd}(p_1,Q_1(y))=1$.
By Remark~3 of Theorem~\ref{RS}, there exist $v\in R^*$ and $Q_5(y)\in R[y]$ such that $\s(x)=vp_3^{-1}(Q_5(p_3x+Q_3(y)-y))$ and $Q_5(Q_3(y))=y$ mod $p_3$.
We have:\\
$(*)\,\,F(x,y)=\s\tau(y)=-vp_3^{-1}p_4y+(vp_3^{-1}p_4Q_5+Q_4)(p_3x+Q_3(y))$.\\
Using $(*)$, we prove that $(p_3,-vp_3^{-1}p_4,Q_4,vp_3^{-1}p_4Q_5+Q_4)\in L_2^{red}(R)$ and this quadruplet is equivalent to $(p_1,p_2,Q_1,Q_2)$.
Uniqueness in Proposition~\ref{P0} gives $p_1=p_3$ and $Q_1(y)=Q_3(y)$ and then $\phi_{p_1}(Q_1(y))\in{\rm VA}_1(R/p_1R)$ by Theorem~\ref{RS}.
Conversely, we assume that $\phi_{p_1}(Q_1(y))\in{\rm VA}_1(R/p_1R)$. Let $Q_5(y)\in R[y]$ be such that $Q_5(Q_1(y))=y$ mod $p_1$.
By Theorem~\ref{RS}, we have: $\s=(p_1^{-1}(Q_5(p_1x+Q_1(y))-y),p_1x+Q_1(y))\in{\rm GA}_2(R)$.
We have: $\s^{-1}(F(x,y))=-p_1p_2x+p_2Q_5(y)+Q_2(y)\in{\cal B}^1(R)\cap{\rm VA}_2(R)$. There exists $\tau\in{\rm GA}_2(R)$ such that
$\tau(y)=\s^{-1}(F(x,y))$ and finally $F=\s\tau(y)$.\\
4) If $F$ is tame $F$ then, by 1), $F$ has length $"1+1"$. If $F$ has a mate of length~$1$ then there exists $\s\in{\rm GA}_2(R)$ such that $\s(y)=G$ and $\s(x)=F$.
If we set: $\tau=\pi$ then $\tau(y)=x\in{\cal B}^1(R)$ and $\s\tau(y)=\s(x)=F$, and $F$ has length $"1+1"$.\\

{\bf Remark.} Let us assume that $R=K[z]$ where $K$ is a field of characteristic~$0$.
Let $F$ be rational length~2 coordinate. Using 1) of Theorem~\ref{P5} one can check if $F$ is tame
coordinate of $K[z][x,y]$. If is not, then Theorem~\ref{UY} implies that $F$ is a wild coordinate of $K[x,y,z]$.

\section{\bf Equivalent polynomials.}

In this section, we assume that $R$ is a UFD.

\begin{de} Let $F,G\in R[x,y]$ we say that $F$ and $G$ are \textit{equivalent} if there exists $\s\in{\rm GA}_2(R)$ such that $\s(F)=G$.
This is of course a equivalent relation.

\end{de}

\begin{theo}\label{Po}
Let $p_1,p_2\in R^{\times}$ be nonzero elements and $Q_1,Q_2\in R[y]$ be polynomials such that
${\rm gcd}(p_1,Q_1(y))={\rm gcd}(p_2,Q_2(y))=1$ and $Q_1(0)=Q_2(0)=0$.\\
1) The polynomials $p_1x+Q_1(y)$ and $p_2x+Q_2(y)$ (in ${\cal B}^1(R)$) are equivalent if and only if there exist $Q_3\in R[y]$ and
$u\in R^*$ such that:
$$(*)\hspace{.3cm}Y=u\,{{\rm gcd}(p_1,p_2)\over p_2}\,\left\{{p_1\over {\rm gcd}(p_1,p_2)}y+Q_3(p_2x+Q_2(y))\right\}\in R[x,y]$$
and $(**)\hspace{.3cm}Q_1(Y)=p_2x+Q_2(y)$ modulo $p_1R[x,y]$.\\
2) If the polynomials $p_1x+Q_1(y)$ and $p_2x+Q_2(y)$ are equivalent then $p_1'x+Q_1(y)$ and $p_2'x+Q_2(y)$ are in ${\rm VA}_2(A)$
where $p_1'=p_1/{\rm gcd}(p_1,p_2)$ and $p_2'=p_2/{\rm gcd}(p_1,p_2)$.
\end{theo}

{\bf Proof.} 1) We assume that there exists $\s\in{\rm GA}_2(R)$ such that $\s(p_1x+Q_1(y))=p_2x+Q_2(y)$. Let $\tau_1=(p_1x+Q_1(y),y)$
and $\tau_2=(p_2x+Q_2(y),y)$ be two triangular automorphisms of $GA_2({\rm qt}(A))$. We have: $\s\tau_1(x)=\tau_2(x)$ {\it i. e.}
$\s\tau_1\pi(y)=\tau_2\pi(y)$ (recall that $\pi=(y,x)$). By Corollary~2, there exist $p_3\in{\rm qt}(R)^*$ and $Q_4\in {\rm qt}(R)[y]$ such that
$\pi\tau_2^{-1}\s\tau_1\pi=\tau_3=(p_3x+Q_4(y),y)$. We set $u=p_2p_3p_1^{-1}$. Since $\s=\tau_2\pi\tau_3\pi\tau_1^{-1}\in{\rm GA}_2(R)$,
we have: $u=\det(J\s)\in R^*$.\\
We set: $Y=\s(y)=\tau_2\pi\tau_3(x)=up_1p_2^{-1}y+Q_4(p_2x+Q_2(y))\in R[x,y]$. Since ${\rm gcd}(p_2,Q_2(y))=1$,
this implies that there exists $Q_3\in R[y]$ such that $Q_4(y)=u\,{\rm gcd}(p_1,p_2)p_2^{-1}Q_3(y)$ and $(*)$ follows.\\
Finally, $\s(x)=\tau_2\pi\tau_3\pi(p_1^{-1}(x-Q_1(y))=p_1^{-1}(\tau_2\pi\tau_3(y)-Q_1(\s(y)))$\\
\hspace*{2.2cm}$=p_1^{-1}(p_2x+Q_2(y)-Q_1(Y))\in R[x,y]$ and we obtain $(**)$.\\
Conversely, if we have $(*)$ and $(**)$, we define an endomorphism $\s$ of $R[x,y]$ by $\s(y)=Y$ and $\s(x)=p_1^{-1}(p_2x+Q_2(y)-Q_1(Y))$.
We can check easily that $\s\in {\rm GA}_2({\rm qt}(R))$ and $\det(J\s)=u\in R^*$. Lemma~\ref{vde} implies that $\s\in{\rm GA}_2(R)$ and a
straight forward computation shows that $\s(p_1x+Q_1(y))=p_2x+Q_2(y)$.\\
2) The assumption $Y={u\over p_2'}\{p_1'y+Q_3(p_2x+Q_2(y))\}\in R[x,y]$ is equivalent to $p_1'y+Q_3(Q_2(y))\in p_2' R[y]$.
Since $p_1'$ is invertible modulo $p_2'$ there exists such $Q_3$ if and only if $Q_2$ is invertible (for composition) modulo $p_2'$.
By Theorem~\ref{RS}, this is equivalent to $p_2'x+Q_2(y)\in{\rm VA}_2(R)$. By symmetry, we have also $p_1'x+Q_1(y)\in{\rm VA}_2(R)$.\\

{\bf Remark\ } Let $\s$ be the automorphism in Theorem~\ref{Po}. Then $Y=\s(y)$ is a rational length~2 coordinate.

\begin{coro}\label{P1}
Let $p_1,p_2\in R^{\times}$ be such that ${\rm gcd}(p_1,p_2)=1$
and $Q_1,Q_2\in R[y]$ be polynomials such that ${\rm gcd}(p,Q_1(y))={\rm gcd}(p,Q_2(y))=1$ and
$Q_1(0)=Q_2(0)=0$. Then $p_1x+Q_1(y)$ and $p_2x+Q_2(y)$ are equivalent if and only both are in ${\rm VA}_2(R)$.
\end{coro}

{\bf Proof.} If $p_1x+Q_1(y)$ and $p_2x+Q_2(y)$ are in ${\rm VA}_2(R)$, both are equivalent to $x$, hence are equivalent.
The converse follows from 2) of Theorem~\ref{Po} (since $p_1=p_1'$ and $p_2=p_2'$).

\begin{coro}\label{P2}
Let $p\in R^{\times}$ be a nonzero element and $Q_1,Q_2\in R[y]$ be polynomials such that ${\rm gcd}(p,Q_1(y))={\rm gcd}(p,Q_2(y))=1$ and
$Q_1(0)=Q_2(0)=0$. There exists $\s\in{\rm GA}_2(R)$ such that $\s(px+Q_1(y))=px+Q_2(y)$ if and only if there exist $Q_3\in R[y]$ and
$u\in R^*$ such that:
$$Q_1(u\{y+Q_3(Q_2(y))\})=Q_2(y)\ {\rm mod}\ pR[y].$$
\end{coro}

{\bf Proof.} Take $p_1=p_2=p$ in 1) of Theorem~\ref{Po}.\\

{\bf Example (Poloni).} We consider, in the ring $R=\C[z]$, the element $p=z^2$, and the polynomials $Q_1(y)=-y^2-zq_1(y)$ and $Q_2(y)=-y^2-zq_2(y)$
where $q_1,q_2\in\C[y]$ are such $q_1(0)=q_2(0)=0$. We have the following characterization:\\
There exist $\s\in{\rm GA}_2(R)$ such that $\s(px+Q_1(y))=px+Q_2(y)$ if and only if $q_2(y)+q_2(-y)=q_1(y)+q_1(-y)$.\\
If we compose the automorphism $\s$ with $(x,y,az)\in{\rm GA}_3(\C)$ where $a\in\C^*$, we obtain the if part of Theorem~4.2.28
p. 96 in [19].\\
In fact, using Corollary~\ref{P2}, there exists $\s\in{\rm GA}_2(R)$ such that $\s(px+Q_1(y))=px+Q_2(y)$  if and only if there
exist $u\in R^*$ and $Q_3\in R[y]$ such that:\\
$(\dag)\,u^2(y+Q_3(-y^2-zq_2(y)))^2+zq_1(u(y+Q_3(-y^2)))=y^2+zq_2(y)\,{\rm mod}\,z^2.$\\
Looking at equation $(\dag)$ modulo $z$, we have:\\
$u^2(y^2+2yQ_3(-y^2)+Q_3(-y^2)^2)=y^2\,\,{\rm mod}\,\,z$, and we deduce: $u^2=1$ and $Q_3(y)=0\,\,{\rm mod}\,\,z$.
We can write $Q_3(y)=zQ_4(y)$ with $Q_4\in\C[z][y]$. Now, $(\dag)$ is equivalent to $2yQ_4(-y^2)=q_2(y)-q_1(uy)$.
There exists such a $Q_4$ if and only if $q_2(y)-q_1(uy)$ is an odd polynomial. We conclude by observing that there exist
$u\in\{-1,1\}$ such that $q_2(y)-q_1(uy)=-q_2(-y)+q_1(-uy)$ if and only if $q_1(y)+q_1(-y)=q_2(y)+q_2(-y)$.

\section{\bf Co-tame automorphisms.}
In this section $K$ is a field of characteristic $0$.\\

We denote by $G={\rm GA}_3(K)$ the group of all automorphisms
of the $K$-algebra $K[x,y,z]$, $A={\rm Aff}_3(K)$ the affine automorphisms sub-group, $B={\rm BA}_3(K)$ the triangular automorphisms sub-group
and $T=<A,B>_G={\rm TA}_3(K)$ the tame automorphisms sub-group.

\begin{de}
Let $\s\in G$ we say that $\s$ is {\it co-tame} if $T\subset <A,\s>_G$. In other words, $\s$ is co-tame, if every tame automorphism
is in the sub-group generated by $\s$ and all affine automorphisms.
\end{de}

{\bf Remark\ } Let $\s\in G$ be an automorphism.\\
\textbf{1)} If $\s$ is tame, we have $T\supset <A,\s>_G$ (this is the origin of our terminology).\\
\textbf{2)} If $\s$ is both tame and co-tame, we have $T=<A,\s>_G$.\\
\textbf{3)} If $\s$ is affine then $\s$ is tame but is not co-tame.\\
\textbf{4)} If $\s$ is co-tame then all automorphisms in $A\s A$ are also co-tame. "To be tame" and "to be co-tame"
are properties of the orbits $A\s A$ ($\s\in G$).\\
\textbf{5)} Let $\alpha\in A$ be an affine automorphism. If $\s\alpha\s^{-1}$ is co-tame then $\s$ is co-tame.
We often use this with $\alpha=t=(x+1,y,z)$.\\

Derksen first proved the existence of a co-tame automorphism (see [12] and Lemma~\ref{Der}).
Bodnarchuk proved that a large class of tame automorphisms are co-tame (see [2] and Theorem~\ref{Bpar}).
They both work in dimension $n\ge 3$ but here we focus in dimension $3$.\\
We denote by $P=\{\s\in G\,;\,\s(y),\s(z)\in K[y,z]\}\subset T$
the set of \textit{parabolic} automorphisms. The automorphisms in $BAB$ (resp. $PAP$) are called bi-triangular
(resp. bi-parabolic).

\begin{lem}[Derksen, 1997]\label{Der}
The automorphism $(x+y^2,y,z)$ is co-tame.
\end{lem}

\begin{lem}[Bodnarchuk, 2004]\label{Btri}
If $\s\in (B\cup BAB)\pri A$ then $\s$ is co-tame.
\end{lem}

\begin{theo}[Bodnarchuk, 2004]\label{Bpar}
If $\s\in (P\cup PAP)\pri A$ then $\s$ is co-tame.
\end{theo}

{\bf Remark\ } Bodnarchuk considers only tame automorphisms. He asks the following question:
Is all non affine tame automorphisms are
co-tame? In other words, is $A$ is a maximal sub-group of $T$? This question is still open.\\

Theorem~\ref{CT1} and Theorem~\ref{CT2} give a lot of examples of non tame automorphisms which are co-tame.

\begin{theo}\label{CT1} Let $\s\in G\pri A$ be a non affine automorphism.
We assume that $\s(z)=z$ and $\s(y)\in{\cal R}_1(K[z])$.
then $\s$ is co-tame.
\end{theo}

{\bf Proof} Let $p_1\in K[z]^{\times}$ and $Q_1(y)\in K[z][y]$ be such that $\s(y)=p_1x+Q_1(y)$.
Since $\s(z)=z$ we can consider $\s$ as an automorphism of $K(z)[x,y]$.
By Theorem~\ref{RS} Remark~3), we have:
$$\s=((up_1)^{-1}(y-Q_2(p_1x+Q_1(y))),p_1x+Q_1(y))$$
where $u\in K^*$ and $Q_2(y)\in K[z][y]$ is such $Q_2(Q_1(y))=y$ mod $p_1K[z][y]$.\\
We can remark that if we consider $\s$ as an automorphism of $K(z)[x,y]$, we have $\s=b_1\pi b_2^{-1}$
with $b_1=(p_1x+Q_1(y),y)$ and $b_2=(up_1x+Q_2(y),y)$.\\
Let $t=(x+1,y,z)\in A$ be the unitary translation on $x$.
An elementary computation gives: $\tau=\s t\s^{-1}=(x+p_1^{-1}(Q_1(y)-Q_1(y+up_1)),y+up_1,z)\in B$.\\
If ${\rm deg}_y(Q_1)\ge 3$ then ${\rm deg}_y(Q_1(y)-Q_1(y+up_1))\ge 2$ and $\tau\in B\pri A$.
By Lemma~\ref{Btri}, $\tau$ and then $\s$ are co-tame.\\
We assume, now, ${\rm deg}_y(Q_1)\le 2$. We write $Q_1(y)=a+by+cy^2$ where $a,b,c\in K[z]$.
We have $p_1^{-1}(Q_1(y)-Q_1(y+up_1))=-u(2cy+b+ucp_1)$.\\
If $c\in K[z]\pri K$ then $\tau\in B\pri A$ and we can conclude as above.\\
We assume, now, $c\in K$. Since $\s(y)\in{\rm VA}_2(K[z])$, Theorem~\ref{RS} implies that
$p_1\in K^*$ or $c=0$. In both cases $\s\in BAB$ and we can conclude using Lemma~\ref{Btri}.\\

{\bf Remark\ } Using Theorem~\ref{SU} and Theorem~\ref{CT1}, we deduce that
 the Nagata automorphism is non tame but co-tame. The sub-group generated by Nagata automorphism and
 affine automorphisms strictly contains the tame automorphisms group! We don't know if this group
 is a proper sub-group of ${\rm GA}_3(K)$.

\begin{theo}\label{CT2} Let $\s\in G\pri A$ be a non affine automorphism.
We assume that $\s(z)=z$ and $\s(y)\in{\cal R}_2(K[z])$.
then $\s$ is co-tame.
\end{theo}

{\bf Proof} There exist $p_1,p_2,p_3\in K(z)^*$ and $Q_1,Q_2,Q_3\in K(z)[y]$ such that
such that $\s=\tau_1\pi\tau_2\pi\tau_3$ where $\tau_i=(p_ix+Q_i(y),y)\in{\rm BA}_2(K(z))$ for
$i\in\{1,2,3\}$ (see Remark 1) after Definition~\ref{Rc}). We prove that $\s$ is co-tame by induction on ${\rm deg}_y(Q_2)$.\\
If ${\rm deg}_y(Q_2)\le 1$ then $\pi\tau_2\pi\in{\rm Aff}_2(K(z))$ and the Bruhat decomposition implies
that $\s(y)\in{\cal R}_1(K[z])$ and $\s$ is co-tame by Theorem~\ref{CT1}.\\
If ${\rm deg}_y(Q_2)\ge 2$ we compute $\tau=\s t\s^{-1}$ where $t=(x+1,y,z)\in A$ is again the unitary translation on $x$.
We have $\tau=\s t\s^{-1}=\tau_1\pi\tau_4\pi\tau_1^{-1}$ where
$\tau_4=\tau_2\pi\tau_3\pi\tau_3^{-1}\pi\tau_2^{-1}=(x+p_2^{-1}(Q_2(y)-Q_2(y+p_3^{-1})),y+p_3^{-1})\in{\rm BA}_2(K(z))$.
Since ${\rm deg}_y(Q_2(y)-Q_2(y+p_3^{-1}))<{\rm deg}_y(Q_2(y))$, by induction $\tau$ and then $\s$ are co-tame.

\section*{References}

\noindent [1]  V. V. Bavula, The inversion formulae for automorphisms of polynomial algebras and rings of differential operators in prime characteristic.
\textit{J. Pure Appl. Algebra}  212  (2008),  no. 10, 2320--2337.

\medskip

\noindent [2] Y. Bodnarchuk, On generators of the tame invertible polynomial maps group, \textit{Internat. J. Algebra Comput.} 15 (2005), no. 5-6, 851–867.

\medskip

\noindent [3] J. Berson, Stable tame coordinates, \textit{J. Pure and Applied Algebra}, 170 (2002), 131-143.

\medskip

\noindent [4] J. Berson, Polynomial coordinates and their behavior in higher dimension, PhD thesis, University of Nijmegen, Nijmegen (2004).

\medskip

\noindent [5] J. Berson, A. van den Essen, S. Maubach, Derivations having divergence zero on $R[X,Y]$, \textit{Israel J. Math.} 124 (2001), 115--124.

\medskip

\noindent [6] S. Bhatwadekar, A. Dutta, Kernel of locally nilpotent $R$-derivations of $R[X,Y]$, \textit{Trans. Amer. Math. Soc.} 349 (1997), no. 8, 3303--3319.

\medskip

\noindent [7] D. Daigle, G. Freudenburg, Locally nilpotent derivations over a UFD and an application to rank two
locally nilpotent derivations of $k[X_1,\cdots,X_n]$, \textit{J. Algebra}  204 (1998), no. 2, 353--371.

\medskip

\noindent [8] V. Drensky and J.-T. Yu, Tame and wild Coordinates of $K[z][x; y]$, \textit{Trans. Amer. Math. Soc.}, 353 (2001), 519-537.

\medskip

\noindent [9] E. Edo, Totally stably tame variables, \textit{J. Algebra} 287 (2005) 15-31.

\medskip

\noindent [10] E. Edo, J.-P. Furter, Some families of polynomial automorphisms, \textit{J. Pure Appl. Algebra} 194 (2004), no. 3, 263--271.

\medskip

\noindent [11] E. Edo and S. V\'en\'ereau, Length 2 variables and transfer, \textit{Annales polonici Math.} 76 (2001), 67-76.

\medskip

\noindent [12] A. van den Essen, Polynomial Automorphisms and the Jacobian Conjecture, \textit{Progress in Mathematics}, vol. 190, Birkh\"auser-Verlag, (2000).

\medskip

\noindent [13] S. Friedland, J. Milnor.  Dynamical properties of plane polynomial automorphisms, \textit{Ergod. Th. Dyn. Syst.} 9 (1989), 67-99.

\medskip

\noindent [14] J-P. Furter, On the variety of automorphisms of the affine plane. \textit{J. Algebra} 195 (1997), no.~2, 604--623.

\medskip

\noindent [15] H. Jung \"Uber ganze birationale Transformationen der Ebene. \textit{J. Reine Angew. Math.} 184 (1942). 161--174.

\medskip

\noindent [16] W. van der Kulk. On polynomial rings in two variables. \textit{Nieuw Arch. Wiskunde} (3) 1 (1953), 33--41.

\medskip

\noindent [17] A. Mikhalev, V. Shpilrain and J.-T. Yu, Combinatorial methods. Free groups, polynomials, and free algebras.
CMS Books in Mathematics 19. Springer-Verlag, New York, 2004. ISBN: 0-387-40562-3.

\medskip

\noindent [18] M. Nagata, On the automorphism group of $k[X,Y]$, \textit{Kyoto Univ. Lectures in Math.} 5 (1972).

\medskip

\noindent [19] P.-M.~Poloni. Sur les plongements des hypersurfaces de Danielewski. Th\`ese de doctorat. Universit\'e de Bourgogne (2008).

\medskip

\noindent [20] P.~Russell. Simple birational extensions of two dimensional affine rational domains. \textit{Compositio Math.} 33 (1976), no. 2, 197--208.

\medskip

\noindent [21] I. Shestakov and U. Umirbaev, The tame and wild automorphisms of polynomial rings in three variables, \textit{J. of the AMS}, 17 (2004), no. 1,
197-227.

\medskip

\noindent [22] U. Umirbaev and  J.-T. Yu, The strong Nagata conjecture. \textit{Proc. Natl. Acad. Sci. USA}, 101 (2004), no. 13.

\medskip
\bigskip

\noindent Author's address:\\

\noindent
Eric Edo\\
edo@univ-nc.nc\\
\noindent
Equipe ERIM,\\
Porte S 20, Nouville Banian,\\
     University of New Caledonia, BP R4,\\
     98 851 Noumea Cedex, New Caledonia.\\

\end{document}